\documentstyle{article}
\begin{document}
\title{A three-parameter deformation of the Weyl-Heisenberg algebra : differential calculus and invariance \thanks{Presented at the 5th Colloquium "Quantum Groups and Integrables Systems", Prague, 20-22 June 1996}}
\author{Mich\`{e}le IRAC-ASTAUD \thanks{e-mail : mici@ccr.jussieu.fr}\\ 
Laboratoire de Physique Th\'{e}orique et Math\'{e}matique\\
Universit\'{e} Paris VII\\2 place Jussieu F-75251 Paris Cedex 05, FRANCE}
\date{}
\maketitle

\begin{abstract}
We define a three-parameter deformation of the Weyl-Heisenberg algebra that generalizes the q-oscillator algebra. By a purely algebraical procedure, we set up on this quantum space two differential calculi that are shown to be invariant on the same quantum group, extended to  a ten-generator Hopf-star-algebra. We prove that, when the values of the parameters are related, the two differential calculi reduce to one that is invariant under two quantum groups.
\end{abstract}

\section{Introduction}
Differential calculus on quantum groups and quantum spaces was developed in many interesting papers \cite{Woro} \cite{wz}. The aim of this talk and of a previous paper \cite{irac} is to set up in a purely algebraical way \cite{wz}\cite{zumino}\cite{boro} a differential calculus on a quadratic deformation of the Weyl-Heisenberg algebra.

The deformation of the Weyl-Heisenberg algebra  considered,  is the free associative algebra generated by  three variables $x^i$ satisfying the  quadratic relations~:

\begin{equation}
R_{xx} :\, 
\left\{
\begin{array}{rcl}
x^1 x^2 - q \, x^2 x^1 -s (x^3)^2 & = & 0 \\
x^1x^3- u \, x^3 x^1 &  = & 0 \\
x^2x^3-u^{-1}\, x^3 x^2 & = & 0
\end{array}
\right.
\label{xx}
\end{equation}
This algebra is by denoted $C<x>/R_{xx}$.

When $q=u^{-2}$,  putting
$x^1=\hat{a}, \quad 
x^2= \hat{a}^{\dagger}$ and $x^3=q^{-\frac{N}{2}}$, the algebra $C<x>/R_{xx}$ is identified to the q-oscillator algebra.

 When $q=1$, $x^3$ is the identity and the q-oscillator algebra restricts to the Weyl-Heisenberg algebra. This algebra is invariant under the seven-dimensional   subgroup, $G$, of 
 $GL(3)$ constituted by the matrices $T$ such as $t^3_1 =t^3_2 =0$ and $t^1_1t^2_2-t^1_2t^2_1=(t^3_3)^2$.	

 When $s=0$, $C<x>/R_{xx}$ is the three-dimensional quantum plane \cite{manin}. 

In \cite{nous}\cite{nous1}, we have proved that a differential calculus on $C<x>/R_{xx}$  invariant on a seven-dimensional quantum group, deformation of $G$, only exists if $q= u^2$ and then is  unique.

 In   section 2, we develop two different differential calculi on $C<x>/R_{xx}$ for arbitrary values of the deformation parameters $q,u$ and $s$ and determine  two sets of  quadratic relations  between variables, differentials and derivatives. In section 3, we show that all these relations are invariant on the same quantum group, that is extended to a ten-generator Hopf-algebra, by adding the inverse of the quantum determinant. When the $x^1$ and $x^2$ are  mutually adjoint and when $x^3$ is self-adjoint, as it is the case for the q-oscillator algebra, the quantum group is endowed with a structure of Hopf-star-algebra. Finally, we discuss 
the particular case where it is a deformation of  $G$, that is when  two elements of the quantum matrix vanish. The two differential calculi (invariant under one quantum group) reduce to one that is invariant under two quantum groups,  deformations of $GL(3)$ and of $G$ \cite{nous}\cite{nous1}.
\section{Differential Calculus} 
Following the general formalism \cite{wz}\cite{zumino}, we add   to the free algebra $C<x>$,
 three generators $\xi^i,\, i=1,2,3$, identified to the one-forms. We define the exterior 
differential operator $d $   in $C<x,\xi>$ such as $ d(x^i)$ is identified to $\xi^i$,  
 $ d^2$ is equal to $0$, $d$ is linear  and  
  satisfies the usual graded  Leibniz rule~:
\begin{equation}
d(fg) = (df)g + (-1)^k f(dg)
\end{equation} 
where $ f,g \in {C<x,\xi>}$ and $f$ is of degree $k$.  
Then, we define the 
partial derivatives such as   $ d \equiv \xi^{i} \; \partial_i$.

The commutation relations $R_{x\xi}$ between the variables and the differentials and those between the partial derivatives $\partial_k$ and the forms  $\xi^l$ are assumed to be quadratic \cite{wz}\cite{zumino}~:
\begin {equation} 
R_{x\xi} \quad: \quad  x^k\xi^l=C^{kl}_{mn} \xi^m x^n 
\label{31}
\end{equation}
and
\begin{equation}
R_{\partial \xi} \quad : \quad
\partial_k \xi^l=K^{lm}_{kn}\xi^n \partial_m
\label{xder} 
\end{equation}
 By applying operator $d$ to (\ref{31})
, we get~:
\begin{equation} 
R_{\xi\xi} \quad :\quad \xi^k \xi^l=-C^{kl}_{mn} \xi^m\xi^n ,\quad \! \! 
\label{32}
\end{equation}
the one-form quantum space corresponds to the C-eigenspace of eigenvalue equal to $-1$. 
Applying the Leibniz rule on $x^k \, f$, $f\in C<x>$, and taking into account
relations (\ref{31}) , we obtain the commutation relations 
between $\partial_i$ and $x^k$ \cite{wz}:
 \begin{equation}
R_{x\partial} \quad :\quad
\partial_lx^k=\delta^k_l+C^{km}_{ln}x^n\partial_m
\label{partialx}
\end{equation}

The combinations of the relations $R_{\partial \xi}, R_{x\partial}$ and $R_{x\xi}$ lead to
\begin{equation}
K  \! \! = \! \!  C^{-1}
\end{equation}
and to a sufficient condition of consistency, the Yang-Baxter equation :
\begin{equation}
(C\otimes 1)(1 \otimes C)(C\otimes 1) = (1 \otimes C)(C\otimes 1)(1 \otimes C)
\label{yb1}
\end{equation}
where 1 is the  identity in $GL(3)$.

From this point, we simplify the matrix $C$ assuming that 
$C^{ij}_{lm}$ is zero if $(i,j)\ne (l, m)$ or $(m,l)$ except $C^{12}_{33}$ and $C^{21}_{33}$.
We multiply  the relations $R_{xx}$  on the left by $\partial_i$. Using $R_{x\partial}$ , we commute $\partial_i$ to the right, and obtain several relations between the elements $C^{ij}_{kl}$. In particular,
\begin{equation}
\begin{array}{llll}
C^{12}_{12}&=q C^{21}_{12}-1,& C^{12}_{21}&=q C^{21}_{21} +q \\
C^{13}_{13}& =u C^{31}_{13} -1,& C^{13}_{31}&= u C^{31}_{31}+u\\
C^{32}_{23}&=u C^{23}_{23} +u,& C^{32}_{32}&=u C^{23}_{32}-1 \\
C^{12}_{33}& = q C^{21}_{33} +s C^{33}_{33} +s&& 
\end {array}
\label{coeff1}
\end{equation}
and 
\begin{equation}
C^{12}_{12} C^{21}_{21} = 0 ,\quad
 C^{13}_{13} C^{31}_{31} = 0 ,\quad
C^{23}_{23} C^{32}_{32} = 0 
\label{coeff2}
\end{equation}
We substitute in the equation (\ref{yb1}) a matrix $C$, the elements of which satisfy (\ref{coeff1}) . This gives $27 \times 27$ new relations between the coefficients $C^{ij}_{kl}$. Solving these relations together with (\ref{coeff2}), we find only two different solutions, namely 
 $C$ is equal to $\Omega$ or to its inverse, with  $\Omega$ defined by~:
\begin{equation}
\Omega =\left(\begin{array}{c c c c c c c c c }
q/u^{2}&0&0&0&0&0&0&0&0\\
0&0&0&q^2/u^2&0&0&0&0&qs/u^2\\
0&0&0&0&0&0&q/u&0&0\\
0&q^{-1}&0&q/u^2-1&0&0&0&0&-s/q\\
0&0&0&0&q/u^2&0&0&0&0\\
0&0&0&0&0&q/u^2-1&0&1/u&0\\
0&0&1/u&0&0&0&q/u^2-1&0&0\\
0&0&0&0&0&q/u&0&0&0\\
0&0&0&0&0&0&0&0&q/u^2\\
\end{array}\right)  
\label{C1}
\end{equation}  
The matrix $\Omega$ and $\Omega^{-1}$ have the same eigenspaces that correspond to the quantum space $C<x>/R_{xx}$, that is our starting point, and to the one-form quantum space defined by~:
\begin{equation} 
R_{\xi\xi} :\,
\left
\{\begin{array}{cccccc}
(\xi^1)^2 & = & 0,\quad & (\xi^2)^2 & = & 0,\quad \\
(\xi^3)^2 & = & 0 \quad &\xi^2 \xi^1& = &-  u^{2}/q^2 \, \xi^1\xi^2 \\  
\xi^1 \xi^3&=&-q/u \, \xi^3 \xi^1, \quad 
& \xi^2 \xi^3 &= & - u/q \, \xi^3 \xi^2 
\end{array}
\right.
\label{xi1}
\end{equation}
When $C$ is equal to $\Omega$ and to its inverse, the eigenspaces of the transpose matrix $(C^{-1})^t$ are the same. The six-dimensional eigenspace
 is identified to the derivative quantum space and is defined by :
\begin{equation} 
R_{\partial \partial} :\,
\left
\{\begin{array}{ll}
 \partial_1 \partial_2 = & u^{2}/q^2 \,\partial_2\partial_1,\\  
\partial_1 \partial_3=&u/q \, \partial_3 \partial_1,\\
\partial_2\partial_3 = &  q/u \, \partial_3 \partial_2 
\end{array}
\right.
\label{partial1}
\end{equation}
The three-dimensional eigenspace corresponds to the covariant differential forms. We denote by $R$ the set of relations (\ref{xx}), (\ref{xi1}) and (\ref{partial1}).

Using (\ref{31}),(\ref{xder}),(\ref{partialx}) with $C=\Omega$ and $C=\Omega^{-1}$, we obtain two  sets of relations  and then two different differential calculi~:

\noindent The commutation relations (\ref{31}) between the variables and the one-forms read

\noindent $\bullet$  when $C=\Omega \, ,\quad R^\Omega_{x\xi}\, :$
\begin{equation}
\begin{array}{llll}
 x^i \xi^i&= q/u^2\,\xi^i x^i,\quad i=1,2,3,&x^1 \xi^3&=q/ u \,\xi^3 x^1,\\
x^1 \xi^2&= q^2/u^2 \,\xi^2 x^1+ qs/u^2 \, \xi^3 x^3,&
x^3 \xi^2&= q/u \, \xi^2 x^3,\\
x^2 \xi^3&= (q/u^2-1) \, \xi^2 x^3+ 1/u \, \xi^3 x^2,&x^3 \xi^1&=( q/u^2-1) \, \xi^3 x^1+ 1/u \, \xi^1 x^3\\
x^2 \xi^1&= 1/q \, \xi^1 x^2+(q/u^2-1) \, \xi^2 x^1- &s/q \,& \xi^3 x^3,\\
\end{array}
\end{equation}
\noindent $\bullet$  when $C=\Omega^{-1}\, ,\quad R^{\Omega^{-1}}_{x\xi}\, :$
\begin{equation}
\begin{array}{llll}
 x^i \xi^i&=  u^2/q \,\xi^i x^i,\quad i=1,2,3,&x^1 \xi^3&= ( u^2/q -1) \, \xi^1x^3+ u \, \xi^3 x^1,\\
x^3 \xi^1&= u/q \, \xi^1x^3, &x^2\xi^1&= u^2/q^2 \, \xi^1x^2- su^2/q^2 \,\xi^3x^3,\\
x^2 \xi^3&= u/q \, \xi^3 x^2 ,&x^3 \xi^2&= ( u^2/q -1) \,\xi^3 x^2+ u  \,\xi^2 x^3,\\
x^1 \xi^2&= ( u^2/q -1) \,\xi^1 x^2 + q  \,\xi^2 x^1 + &s \,\xi^3x^3.&
\end{array}
\end{equation}
The commutation relations (\ref{xder}) between the derivatives and the differentials

\noindent $\bullet$  when $C=\Omega \, ,\quad R_{\partial \xi}^{\Omega} \, :$ 
\begin{equation}
\begin{array}{llll}
\partial_3 \xi^3&= (u^2/q -1)\,\xi^2\partial_2+ u^2/q \, \xi^3\partial_2
&\partial_1 \xi^2&=  u^2/q^2 \,\xi^2\partial_1,\\
\partial_1 \xi^3&=  u/q \,\xi^3\partial_1,
&\partial_2 \xi^1&= q \, \xi^1\partial_2,\\
\partial_3 \xi^2&=  u/q \,\xi^2\partial_3- su^2/q^2 \,\xi^3\partial_1,
&\partial_2 \xi^3&= u \, \xi^3\partial_2,\\
\partial_3 \xi^1&= u \, \xi^1\partial_3+ s \, \xi^3\partial_2,
&\partial_2 \xi^2&=  u^2/q \, \xi^2\partial_2,\\
\partial_1 \xi^1&=  u^2/q \, \xi^1\partial_1+(u^2/q -1)\,\ \xi^3\partial_3+&(u^2/q -1)&\,\ \xi^2\partial_2,
\end{array}
\end{equation}
\noindent $\bullet$  when $C=\Omega^{-1}\, ,\quad R_{\partial \xi}^{\Omega^{-1}}\, :$ 
\begin{equation}
\begin{array}{llll}
 \partial_1 \xi^1&= q/u^2\,\xi^1\partial_1,&\partial_3 \xi^2&= 1/u \, \xi^2\partial_3- s/q \, \xi^3\partial_1,\\
\partial_1 \xi^3&= 1/u \, \xi^3\partial_1,
&\partial_2 \xi^1&= q^2/u^2\,\xi^1\partial_2,\\
\partial_2 \xi^3&= q/u \,\xi^3\partial_2,&\partial_3 \xi^1&= q/u \, \,\xi^1\partial_3+ sq/u^2\, \,\xi^3\partial_2,\\
\partial_1 \xi^2&= 1/q \, \xi^2\partial_1,
&\partial_3 \xi^3&= (q/u^2-1) \,\xi^1\partial_1+ q/u^2\, \xi^3\partial_3,\\
\partial_2 \xi^2&= (q/u^2-1) \,\xi^1\partial_1+&(q/u^2-1) \,&\xi^3\partial_3+ q/u^2\,\xi^2\partial_2
\end{array}
\end{equation}
The commutation relations (\ref{partialx}) between the derivatives and the variables

\noindent $\bullet$  when $C=\Omega \, ,\quad R_{x \partial}^{\Omega} \, :$ 
\begin{equation}
\begin{array}{llll}
  \partial_1 x^1&= 1+  q/u^2\, x^1\partial_1,&\partial_2 x^3&=  q/u \,x^3\partial_2,\\
\partial_3 x^3&= 1+  q/u^2\,x^3\partial_3+ (q/u^2-1) \, x^1\partial_1,&\partial_1 x^2&= 1/q \, x^2\partial_1,\\
\partial_3 x^1&=  q/u \,x^1\partial_3+  qs/u^2\, x^3\partial_2,
&\partial_2 x^1&=  q^2/u^2\,x^1\partial_2,\\
\partial_3 x^2&= 1/u \, x^2\partial_3-s/q \,x^3\partial_1,&\partial_1 x^3&= 1/u \, x^3\partial_1,\\
\partial_2 x^2&= 1+ q/u^2\,x^2\partial_2+  (q/u^2-1) &\,x^1\partial_1+& (q/u^2-1) \, x^3\partial_3.\\
\end{array}
\end{equation}
\noindent $\bullet$  when $C=\Omega^{-1} \, ,\quad R_{x \partial}^{\Omega^{-1}} \, :$ 
\begin{equation}
\begin{array}{llll}
\partial_2 x^2&=1+ u^2/q \,x^2\partial_2, 
&\partial_1 x^3&= u/q \,x^3\partial_1,\\
\partial_3 x^3&= 1+ u^2/q \,x^3\partial_3+ (u^2/q -1) \,x^2\partial_2,&
\partial_2 x^1&= q \, x^1\partial_2,\\
\partial_3 x^2&= u/q \,x^2\partial_3- su^2/q \,x^3\partial_1,&\partial_1 x^2&= u^2/q^2 \, x^2\partial_1,\\
\partial_3 x^1&= u \,x^1\partial_3 +s \, x^3 \partial_2,
&\partial_2 x^3&= u \, x^3\partial_2,\\
\partial_1 x^1&= 1+ u^2/q \, x^1\partial_1+(u^2/q -1) \,x^2\partial_2+&(u^2/q -1) \,& x^3\partial_3.\\
\end{array}
\end{equation}
The two sets of relations $(R, R^{\Omega})$ and $(R, R^{\Omega^{-1}})$ define two quadratic algebras, $C<x,\xi,\partial>/R \cup R^{\Omega}$ and 
$C<x,\xi,\partial>/R \cup R^{\Omega^{-1}}$. In the following section, we investigate their invariance.

It is to be noted that all the  construction of the differential calculus is performed without using  the B-matrix associated with the variables \cite{wz} \cite{zumino}, and is the result solely  of the relations $R_{xx}, R_{x \xi}$ and $R_{\partial \xi}$ . Moreover, as a consequence of the construction, $B$ is found to be equal to $C$. 
\section{Quantum Group and Invariance}
The quantum matrix $T$ with nine non commuting elements defines a homomorphism on $C<x,\xi,\partial>$ \cite{rtf}. The  variables $x$ and  the differentials $\xi$ are transformed by $T$ and the derivatives $\partial$ are transformed by $(T^{-1})^ t$.

\noindent When the matrix $T$ satisfies 
\begin{equation}
{R}^{ji}_{kl} \, t^k_m t^l_n = t^j_l t^i_k \, {R}^{lk}_{mn}
\end{equation}
with $R=\Omega$ (resp. $R=\Omega^{-1}$), the  relations $R \cup R^{\Omega}$ (resp. $R \cup R^{\Omega^{-1}}$) are invariant, and therefore  this homomorphism maps $C<x,\xi,\partial>/R \cup R^{\Omega}$ (resp. $C<x,\xi,\partial>/R \cup R^{\Omega^{-1}}$) in itself.
  It is easy to see that $R=\Omega$ and $R=\Omega^{-1}$ define the same quantum matrix $T$, the elements of which satisfy the following commutation relations $R_{tt}$ :
 
\begin{equation}
\begin{array}{llll}
t^1_2t^1_1&=q^2/u^2\, t^1_1t^1_2,&t^2_2t^1_1&=t^1_1t^2_2-(u^2-q)/q^2\, t^1_2t^2_1- qs /u^2\, t^3_1t^3_2,\\
t^1_3t^1_2&= u /q \,t^1_2t^1_3,&t^2_1t^1_1&=1 /q \,t^1_1t^2_1 - s/q \,(t^3_1)^2,\\
t^1_3t^1_1&= q/u \,t^1_1t^1_3,&t^2_3t^1_1&= u/q \,t^1_1t^2_3-(u^2-q)/q^2 \,t^1_3t^2_1- s/q \,t^3_3t^3_1,\\
t^3_2t^2_2&=u \,t^2_2t^3_2,&t^3_3t^1_1&=t^1_1t^3_3- (u^2-q)/(uq)\, t^1_3t^3_1,\\
t^3_1t^1_1&= 1/u\, t^1_1t^3_1,&t^3_2t^1_1&= q/u \,t^1_1t^3_2- (u^2-q)/u\, t^1_2t^3_1,\\
t^3_3t^1_2&=1/q \,t^1_2t^3_3,&t^2_2t^1_2&= 1/q \,t^1_2t^2_2- s/q \,(t^3_2)^2,\\
t^2_3t^2_2&= u/q \,t^2_2t^2_3,&t^2_3t^1_2&= u/q^2\, t^1_2t^2_3- s/q \,t^3_3t^3_2,\\
t^3_1t^1_2&= u/q^2\, t^1_2t^3_1,&t^2_1t^1_2&= u^2/q^3 \,t^1_2t^2_1- s/q\, t^3_1t^3_2,\\
t^3_2t^2_1&= q^2/u \,t^2_1t^3_2,&t^3_2t^1_3&=t^1_3t^3_2- (u^2-q)/(uq)\, t^1_2t^3_3,\\
t^3_2t^1_2&=1/u \,t^1_2t^3_2, &t^2_3t^1_3&= 1/q\, t^1_3t^2_3- s/q \,(t^3_3)^2+ s/q \,t^1_1t^2_2- su^2/q^3\, t^1_2t^2_1,\\
t^3_1t^1_3&= 1/q\, t^1_3t^3_1,&t^2_2t^1_3&= u/q\,t^1_3t^2_2- (u^2-q)/q^2\, t^1_2t^2_3- s/u \,t^3_3t^3_2,\\
t^2_1t^2_2&= u^2/q^2\, t^2_2t^2_1,&t^3_3t^2_2&=t^2_2t^3_3+ (u^2-q)/u \,t^2_3t^3_2\\
t^2_3t^2_1&= q/u \,t^2_1t^2_3,&t^3_3t^2_3&=u \,t^2_3t^3_3+ sq/u \,t^2_1t^3_2-su\,t^2_2t^3_1,\\
t^3_2t^3_1&= q^2/u^2\, t^3_1t^3_2,&t^3_1t^2_2&= u/q \,t^2_2t^3_1 + (u^2-q)/u \,t^2_1t^3_2,\\
t^3_2t^3_3&= q/u \,t^3_3t^3_2,&t^3_1t^2_3&=t^2_3t^3_1+ (u^2-q)/u \,t^2_1t^3_3,\\
t^3_1t^2_1&=u \,t^2_1t^3_1,& t^2_1t^1_3&= u/q^2 \,t^1_3t^2_1 - su/q^2 \,t^3_3t^3_1\\
t^3_1t^3_3&= u/q \,t^3_3t^3_1,&t^3_3t^1_3&=1/u \,t^1_3t^3_3 +s/u\, t^1_1t^3_2-su/q^2 \,t^1_2t^3_1,\\
t^3_3t^2_1&=q \,t^2_1t^3_3,&t^3_2t^2_3&=q \,t^2_3t^3_2\\
\label{RT}
\end{array}
\end{equation}
Any elements  of $C<t>/R_{tt}$ can be written as a sum of ordered monomials $(t^1_1)^{k_1} (t^1_2)^{k_2}(t^1_3)^{k_3}(t^2_2)^{k_4}(t^2_1)^{k_5}(t^2_3)^{k_6}(t^3_3)^{k_7}(t^3_1)^{k_8}(t^3_2)^{k_9}$ by using the relations $R_{tt}$.
The inverse $T^{-1}$, of $T$ is equal to :
\begin{equation}
\left(\begin{array}{c c c}
t^2_2t^3_3-q/u \,t^2_3t^3_2&-q^2/u^2\, t^1_2t^3_3 + q^3/u^3\, t^1_3t^3_2&t^1_2t^2_3-q/u \,t^1_3t^2_2\\
-u^2/q^2\, t^2_1t^3_3 +u^3/q^3 \,t^2_3t^3_1&t^1_1t^3_3-u/q \,t^1_3t^3_1&-u^2/q^2 \,t^1_1t^2_3+ u^3/q^3 \,t^1_3t^2_1\\
t^2_1t^3_2 -u^2/q^2 \,t^2_2t^3_1&-q^2/u^2 \,t^1_1t^3_2 +t^1_2t^3_1&t^1_1t^2_2 -u^2/q^2\, t^1_2t^2_1\\
\end{array}\right) D^{-1} 
\label{t1}
\end{equation}
where $D$ is the quantum determinant  equal to 
$$D =t^1_1t^2_2t^3_3+t^1_3t^2_1t^3_2+u^3/q^3 \,t^1_2t^2_3t^3_1- q/u \,t^1_1t^2_3t^3_2-u^2/q^2 \,t^1_2t^2_1t^3_3 -u^2/q^2\, t^1_3t^2_2t^3_1.$$ 
We  calculate and verify that $D$ is not a central element of $C<t>$. The elements of $T^{-1}$ belong to the algebra generated by  $D^{-1}$ and by the  $t^i_j$. We deduce the commutation relations $R_{t D^{-1}} $ of $D^{-1}$ with all the $t^i_j$  from those of $D$~:
\begin{equation}
\begin{array}{cccccc}
t^1_1 D^{-1}&=D^{-1} t^1_1,& t^1_2 D^{-1}&=u^2/q^4 t^1_2\,D^{-1},& t^1_3 D^{-1}&= u/q^2 \,D^{-1}t^1_3,\\
t^2_2 D^{-1}&=D^{-1} t^2_2,& t^2_1 D^{-1}&= q^2 \,D^{-1}t^2_1,& t^2_3 D^{-1}&= u/q^2 \,t^2_3 D^{-1},\\
t^3_1 D^{-1} &= q^2/u \,D^{-1} t^3_1,& t^3_2 D^{-1}&=u/q^2\, t^3_2 D^{-1},& t^3_3 D^{-1}&= D^{-1}t^3_3.\\
\end{array}
\label{RD}
\end{equation}
The quotient algebra $C< t, D^{-1}>/R_{tt} \cup R_{t D^{-1}}$ is a Hopf algebra with the co-product $\Delta$, the co-unit $\epsilon$ and antipode $S$  defined by~:  
\begin{eqnarray}
\Delta(T) \equiv T \otimes T, \quad 
\Delta(D^{-1}) \equiv D^{-1} \otimes D^{-1}
\label{26} \\
\epsilon (T,D^{-1}) \equiv (I,1), \quad 
S(T)\equiv T^{-1}, \quad   S(D) \equiv D^{-1}
\label{28}
\end{eqnarray}

When $x^1$ and $x^2$ are mutually adjoint and when $x^3$ is self-adjoint (for instance, in the case of the q-oscillator algebra ), the relations $R_{xx}$ are unchanged if the parameters are real. The action of $T$ respects this property if the quantum group is equipped with a star-operation that is an antihomomorphism such as 
\begin{equation}
\begin{array}{llllll}
&((t^i_j)^*)^*&=t^i_j,&\quad &\forall i,j&\\
(t^2_2)^*&=t^1_1,&\quad (t^2_1)^*&=t^1_2,&\quad
 (t^2_3)^*&=t^1_3,\\
&( t^3_1)^*&=t^3_2,&\quad (t^3_3)^*&=t^3_3&.
\end{array}
\end{equation}
These relations are consistent with (\ref{RT}) and (\ref{RD}) and the quantum group $C<t,D^{-1}>/R_{tt}\cup R_{tD^{-1}}$ acquires the structure of a Hopf-star-algebra.

Let us give two interesting particular cases of $C<t,D^{-1}>/R_{tt}\cup R_{t D^{-1}}$~:

$\bullet$ 
When $s=0$, the variable quantum space is the three-dimensional quantum plane, the resulting quantum group  corresponds to an original deformation of $GL(3)$.

$\bullet$ 
When $t^3_1$ and $t^3_2$ vanish, $C<t,D^{-1}>/R_{tt}\cup R_{t D^{-1}}$ becomes a deformation $G_{q s}$ of the subgroup $G$ of $GL(3)$.
      Two of the relations $R_{tt}$ (\ref{RT}) give  : 
$$(u^2 -q)\, t^1_2 t^3_3=0,$$
 and 
 $$(u^2 -q)\, t^2_1 t^3_3=0.$$ 
implying that  $q$ is equal to $u^2$ if the algebra has no zero divisors. When $q=u^2$,
 the matrix $\Omega$ being equal to its inverse,  the two differential calculi on $C<x>/R_{xx}$ reduce to one. This differential calculus was previously obtained by a completely different method, implying the uniqueness of the result \cite{nous}\cite{nous1}.  All the commutation relations are invariant on the ten-generator quantum group $C<t,D^{-1}>/R_{tt}\cup R_{tD^{-1}}$ and on a quantum group $G_{qs}$ deduced from the previous one  by putting $t^3_1=t^3_2=0.$

In the case where $q=u^2$ we would point out that, if we add the generator $(t^3_3)^{-1}$ to the quantum group $G_{q s}$, the $T$ matrix can be written on the form $T^\prime \times t^3_3$ 
with $t^{\prime i}_j = t^i_j (t^3_3 )^{-1}$. All the elements $t^{\prime i}_j$ commute two by two. Nevertheless, they cannot be identified to C-numbers (and then $T^\prime$ to a matrix belonging to the initial subgroup $G$), because if they were C-numbers all the elements $t^i_j$ would be proportional to $t^3_3$ and this is impossible due to their commutation relations resulting from (\ref{RT}).
\section{Conclusion}
Two differential calculi can be associated with the three-parameter oscillator algebra $C<x>/R_{xx}$ and in particular with the q-oscillator algebra. They are invariant on the same quantum group $C<t,D^{-1}>/R_{tt}\cup R_{tD^{-1}}$ that is an original three-parameter deformation of $GL(3)$. When we assume that the variables $x^1$ and $x^2$  are  mutually adjoint and that $x^3$ is self adjoint, a star-operation is defined on the invariance quantum group that  then becomes a Hopf-star-algebra. We have considered the case where two generators, $t^3_1$ and $t^3_2$ are removed from the algebra $C< t, D^{-1}>/R_{tt} \cup R_{t D^{-1}}$, the subalgebra $G_{qs}$ thus obtained exists only if the parameters are related by $q=u^2$ and is a deformation of the invariance subgroup of the Weyl-Heisenberg algebra. The two differential calculi then reduce to one invariant under two quantum groups, $C< t, D^{-1}>/R_{tt} \cup R_{t D^{-1}}$ and $G_{qs}$.

We would like to thank J.Bertrand for  many interesting discussions.

 \end{document}